\documentclass[12pt]{article}
\usepackage{latexsym}
\usepackage{amsmath,amssymb,amsthm}

\usepackage{amsmath,amsthm,latexsym,amssymb,pgf,makeidx,authblk,boxedminipage,tikz}

\usetikzlibrary{arrows,decorations.pathmorphing,backgrounds,positioning,fit,petri,calc}

\begin{document}

% \author{Peter Dankelmann}

% \email{pdankelmann@uj.ac.za}
% \subjclass[2010]{Primary 05C12; secondary 05C20, 05C62}

%\begin{titlepage}
\title{The Steiner $k$-Wiener index of graphs with given minimum degree} 
\author{Peter Dankelmann (University of Johannesburg)}
\maketitle

\begin{abstract}
Let $G$ be a connected graph. The Steiner distance $d(S)$ of a set $S$ of vertices
is the minimum size of a connected subgraph of $G$ containing all vertices
of $S$. For $k\in \mathbb{N}$, the Steiner $k$-Wiener index $SW_k(G)$ is defined as 
$\sum_S d(S)$, where the sum is over all $k$-element subsets of the 
vertex set of $G$. The average Steiner $k$-distance $\mu_k(G)$ of $G$ is defined as
$\binom{n}{k}^{-1} SW_k(G)$. 

In this paper we prove upper bounds on the Steiner Wiener index and the average 
Steiner distance of graphs with given order $n$ and minimum degree $\delta$. Specifically 
we show that 
$SW_k(G) \leq \frac{k-1}{k+1}\frac{3n}{\delta+1} \binom{n}{k} + O(n^{k})$, 
and that $\mu_k(G) \leq \frac{k-1}{k+1}\frac{3n}{\delta+1} + O(1)$. 
We improve this bound for triangle-free graphs to
$SW_k(G) \leq \frac{k-1}{k+1}\frac{2n}{\delta} \binom{n}{k} + O(n^{k})$, 
and $\mu_k(G) \leq \frac{k-1}{k+1}\frac{2n}{\delta} + O(1)$. 
All bounds are best possible. 
\end{abstract}
MSC-class: 05C12 (primary) 92E10 (secondary)
%\end{center}

%\end{titlepage}

\newtheorem{defi}{Definition}
\newtheorem{theo}{Theorem}
\newtheorem{coro}{Corollary}
\newtheorem{prop}{Proposition}
\newtheorem{rem}{Rem:}
\newtheorem{la}{Lemma}
\newtheorem{exa}{Example}
\newtheorem{conj}{Conjecture}
\newtheorem{problem}{Problem}

%\begin{center}
%\Large{The Steiner $k$-Wiener index of graphs with given minimum degree} \\[2mm]
%\large {Peter Dankelmann} \\
%University of Johannesburg
%\end{center}

% \normalsize

% Keywords: Steiner Wiener index; weighted Steiner Wiener index, 
% average Steiner distance; Wiener index; transmission, average distance; 
% mean distance, Steiner distance.

\section{Introduction}

The Wiener index $W(G)$ of a connected graph $G$ is defined as the sum of the
distances between all unordered pairs of vertices, i.e.,
\[ W(G) = \sum_{\{u,v\} \subseteq V(G)} d(u,v), \]
where $V(G)$ is the vertex set of $G$, and $d(u,v)$ is the usual distance,
i.e., the minimum length of a path from $u$ to $v$. First studied by
the chemist Wiener \cite{Wie1947} as an indicator for the boiling point
of certain alkanes, the Wiener index has become one of the most important
topological indices in chemical graph theory. For its many applications
see, for example, the survey \cite{Rou2002}. 
Since its inception, the Wiener index has attracted much interest in the graph
theory literature under different names, such as transmission, defined as the 
sum of the distances between all ordered pairs of vertices, and total distance. 
The Wiener index is
closely related to the average distance $\mu(G)$, also known as mean distance, 
defined as the average of all distances between two vertices of $G$, i.e.,  
\[ \mu(G) = \binom{n}{2}^{-1} \sum_{\{u,v\} \subseteq V(G)} d(u,v),\]  
where $n$ is the order of $G$. Hence $W(G) = \binom{n}{2}\mu(G)$. 

The Steiner distance $d_G(S)$ of a set $S$ of vertices in a connected graph
$G$ is defined as the minimum size of a connected subgraph of $G$
containing all vertices of $S$. This concept 
was introduced by Chartrand, Oellermann, Tian and Zou \cite{ChaOelTiaZou1989}
in order to generalise the notion of distance between two vertices in a graph
to an arbitrary number of vertices. For $k\in \mathbb{N}$, the maximum value 
of $d_G(S)$, taken over all $k$-sets of vertices of $G$,
 is known as the $k$-Steiner diameter or $k$-diameter
of $G$. For results on the Steiner diameter see, for example, 
\cite{AliDanMuk2012, Ali2013, AliDanMuk2014, DanOelSwa1996-2} and the survey
paper \cite{Mao2015}. 

This paper is concerned with the Steiner $k$-Wiener index, which generalises 
the Wiener index by combining the notions of Steiner distance and Wiener index. 
For $k\in \mathbb{N}$, the Steiner $k$-Wiener index $SW_k(G)$ of a
connected graph $G$ is defined as the sum of the Steiner distances of all
$k$-sets of vertices, i.e.,
\[ SW_k(G) = \sum_{S \subseteq V(G), |S|=k} d_G(S). \]
It generalises the Wiener index since clearly $SW_2(G)=W(G)$. The Steiner
$k$-Wiener index was introduced by Li, Mao and Gutman \cite{LiMaoGut2016}.
A closely related graph parameter, the average Steiner $k$-distance of $G$,
denoted by $\mu_k(G)$ and 
defined as the average of the Steiner distances of all $k$-element subsets
of $V(G)$, was introduced in \cite{DanOelSwa1996} and further investigated 
in \cite{DanOelSwa1997}. In the same way in which the Steiner $k$-Wiener
index generalises the Wiender index, the average $k$-distance generalises the 
average distance since $\mu_2(G) = \mu(G)$.

Several results on the Wiener index are known to hold also for 
the Steiner $k$-Wiener index. 
The observation that the Wiener index of a graph of order $n$ is at least
$\binom{n}{2}$ easily extends to the Steiner $k$-Wiener index: since the
Steiner distance of a set of $k$ vertices is at least $k-1$, we 
obtain the lower bound $SW_k(G) \geq (k-1)\binom{n}{k}$. 
Doyle and Graver \cite{DoyGra1977}, Entringer, Jackson and Snyder \cite{EntJacSny1976} 
and Lov\'{a}sz \cite{Lov1979} independently observed that 
\begin{equation}   \label{eq:path-maximises-W}
W(G) \leq \frac{n+1}{3} \binom{n}{2}, 
\end{equation}
with equality if and only if $G$ is a path. Dankelmann, Oellermann and Swart
\cite{DanOelSwa1996} and Li, Mao, Gutman \cite{LiMaoGut2016} showed that 
this result extends to the Steiner $k$-Wiener index. 

\begin{theo} \label{theo:bounds-on-Wk-for-trees}
{\rm \cite{DanOelSwa1996}, \cite{LiMaoGut2016}} 
Let $G$ be a graph of order $n$. Then
\[ SW_k(G) \leq \frac{(k-1)(n+1)}{k+1} \binom{n}{k}. \]
Equality holds if $G$ is a path. 
\end{theo}

It was observed by Plesn\'\i k \cite{Ple1984} that among all trees of
given order the star minmimises the Wiener index. It was shown by
Dankelmann, Oellermann and Swart \cite{DanOelSwa1996}, and 
by Li, Mao and Gutman \cite{LiMaoGut2016} that the same statement holds 
also for the Steiner $k$-Wiener index. Plesn\'\i k \cite{Ple1984} also showed
that among all trees of order $n$ and diameter $d$, the Wiener index 
is minimised by the tree obtained from a path on $d+1$ vertices by
attaching $n-d-1$ vertices to a centre vertex of the path. This result
was shown in \cite{LuHuaHouChe2018} to hold for the Steiner $k$-Wiener index.
Nordhaus-Gaddum type results on the Wiener index \cite{ZhaWu2005}
were generalised to the Steiner $k$-Wiener index \cite{MaoWanGutLi2017}, 
and so were results on the inverse Wiener problem, i.e., the question
which numbers are the Wiener index of some graph \cite{Wag2006, WanYu2006, LiMaoGut2018},
as well as results on product graphs (see \cite{YehGut1994, MaoWanGut2016}).

The bound on the Wiener index \eqref{eq:path-maximises-W} has been
improved for graphs with various given properties. For example
for $2$-connected graphs (i.e., connected graphs in which removing
a vertex does not disconnect the graph), Plesn\'\i k \cite{Ple1984}
showed that 
\begin{equation} \label{eq:max-W-of-2-connected}
 W(G) \leq \frac{n}{2} \lfloor \frac{n^2}{4} \rfloor, 
\end{equation}
with equality if and only if $G$ is a cycle. An extension of this
result was given in  \cite{DanOelSwa1996-2}, where it was shown that 
the cycle maximises, among all $2$-connected graphs of given order, the 
Steiner $k$-Wiener index for every integer $k$ with $2 \leq k\leq n$. 
Inequality \eqref{eq:max-W-of-2-connected} was strengthened for graphs
with given order and connectivity \cite{DanMukSwa2009} and edge-connectivity
\cite{DanMukSwa2008, DanMukSwa2008-2}, but no bounds on the 
Steiner $k$-Wiener index of graphs of given connectivity or 
edge-connectivity appear to be known.

The bound in \eqref{eq:path-maximises-W} has also been improved for graphs of given minimum
degree. The computer program GRAFFITI conjectured that the average distance
of a graph of order $n$ and minimum degree $\delta$ is not more than $\frac{n}{\delta}$. 
Kouider and Winkler \cite{KouWin1997} proved the following result which 
is asymptotically sharp. 

\begin{theo} {\rm \cite{KouWin1997}} \label{theo:KouiderWinkler}
Let $G$ be a connected graph of order $n$ and minimum degree $\delta$.
Then 
\[ W(G) \leq (\frac{n}{\delta+1}+2) \binom{n}{2}. \]
\end{theo}

Although Theorem \ref{theo:KouiderWinkler} is asymptotically stronger than 
the above-mentioned GRAFFITI conjecture, it does not actually imply it.
The GRAFFITI conjecture, 
was proved later by Beezer, Riegsecker and Smith \cite{BeeRieSmi2001}. 
Kouider and Winkler's result was extended in two ways by Dankelmann and Entringer \cite{DanEnt2000} by showing that the above bound holds not only for $G$ 
but for some spanning tree of $G$, and further by improving this 
bound by a factor of about $\frac{2}{3}$ for triangle-free graphs. 
It is the aim of this paper to show that these bounds on the 
Wiener index of graphs in terms of order and minimum degree extend to
the Steiner $k$-Wiener index and the average Steiner $k$-distance.

We note that there are several recent bounds on the Wiener index, for example in terms
of diameter \cite{MukVet2014}, radius \cite{DasNad2017}, 
% for Eulerian graphs \cite{GutCruRad2014}, 
and other results \cite{KlaNad2014, KnoLuzSkrGut2014, KnoSkrTep2016} 
on the Wiener index 
which have not (yet) been shown to hold for the Steiner $k$-Wiener index.

\section{Notation}

The notation we use is as follows. By $G$ we always denote a finite, simple, 
connected graph on $n(G)$ vertices with vertex set $V(G)$. For a vertex $v$ of $G$, $N_G(v)$ is the
neighbourhood of $v$, i.e., the set of vertices adjacent to $v$, and $N_G[v]$ 
is the closed neighbourhood of $v$, i.e., the set $N_G(v)\cup\{v\}$. 
For $A \subseteq V(G)$ we define $N[A]=\bigcup_{v\in A} N[v]$. 
The degree ${\rm deg}_G(v)$ of $v$ is the number of vertices in $N_G(v)$,
and the minimum degree $\delta(G)$ of $G$ is the smallest of the degrees of
the vertices of $G$. If $U\subseteq E(G)$, then 
$V(U)$ is the set of vertices of $G$ incident with at least one edge in $U$, and $G[U]$, the subgraph induced
by $U$, is the subgraph whose vertex set is $V(U)$ and whose edge set is $U$. 

The distance between two vertices $u$ and $v$, i.e., the minimum length of a 
$(u,v)$-path, is denoted by $d_G(u,v)$. If $S$ is a non-empty subset of 
the vertex set of $G$, then the Steiner distance of $S$, $d_G(S)$, is the 
minimum size of a connected subgraph of $G$ containing the vertices of $S$. 
The distance between a vertex $v$ and a set $A$ of vertices of $G$ is
defined as $\min_{w \in A} d_G(v,w)$. 
If the graph is understood from the context, then
we sometimes omit the argument or subscript $G$.

By $K_n$ we mean the complete graph on $n$ vertices, and $nK_1$ denotes the
edgeless graph on $n$ vertices. For disjoint graphs
$G_1, G_2,\ldots, G_k$ the sequential sum $G_1+G_2+\cdots+G_k$ is the 
graph obtained from the union of $G_1, G_2, \ldots, G_k$ by joining 
every vertex of $G_i$ to every vertex of $G_{i+1}$ for $i=1,2,\ldots,k-1$. 
The line graph of a graph $G$ is the graph whose vertices are
the edges of $G$, with two vertices of the line graph being
adjacent if the corresponding edges of $G$ have a vertex
in common.  
% If in a sequential sum a pattern
% is repeated $\ell$ times, then we indicate this with
% square brackets and the exponent $\ell$; for example
% $G_1 + [G_2 + G_3]^{\ell} +G_4$ stands for 
% $G_1 + G_2 + G_3 + G_2 + G_3 + \cdots + G_2 + G_3 +G_4$,
% where the pattern $G_2+G_3$ appears $\ell$ times. 

\section{Weighted Steiner $k$-Wiener index}

In this section we introduce the weighted Steiner $k$-Wiener index, which 
generalises the weighted Wiener index. The main result
of this section, Lemma \ref{la:fundamental-lemma}, is a common 
generalisation of Theorem \ref{theo:bounds-on-Wk-for-trees} and 
a bound on the weighted Wiener index given in \cite{DanEnt2000}. 
For the proof of Lemma \ref{la:fundamental-lemma} we require the following
notation. 

\begin{defi}
Given a set $X$ and $c$ a weight function $c: X \rightarrow \mathbb{N}_0$. \\
(a) For $Y \subseteq X$ we define $c(Y)$ as $\sum_{y \in Y} c(y)$. \\
(b) Define $X_c$ to be the set obtained from $X$ by replacing every $x \in X$ with 
$c(x)\geq 1$ by elements $x^1, x^2,\ldots, x^{c(x)}$ and deleting all
$y\in X$ for which $c(y)=0$. If $c(x)>0$ then we refer to
$x$ as the original, and to $x^1, x^2,\ldots, x^{c(x)}$ as the copies
of $x$.  \\
(c) Given a set $Y \subseteq X_c$, the original of $Y$  is the
set $Y^* \subset X$ whose elements are exactly those $x\in X$ of which $Y$ 
contains at least one copy.  \\
(d) If $X$ is the vertex set of a graph $G$, then for $Y \subseteq X_c$
we define $d(Y)$ as $d_G(Y^*)$. 
\end{defi}

It is easy to see that $|X_c| = c(X)$. With the above definition we can
now introduce the weighted Steiner $k$-Wiener index of a graph. 

\begin{defi}
Let $G$ be a connected graph with vertex set $V$ and 
$c: V \rightarrow \mathbb{N}_0$ a weight function. 
The Steiner $k$-Wiener index of $G$ 
with respect to $c$ is defined by
\[ SW_k(G,c) = \sum_{S \subseteq V_c, |S|=k} d_G(S). \]
\end{defi}

  \begin{figure}[h]
  \begin{center}
\begin{tikzpicture}
  [scale=0.7,inner sep=1mm, % this is the node radius
   vertex/.style={circle,thick,draw}, % this defines the default style for the vertex class
   thickedge/.style={line width=2pt}] % this defines the default style for the thick edges
    \begin{scope}[>=triangle 45]
    
     \node[vertex] (a1) at (3.75,0) [fill=white] {};
     \node[vertex] (a2) at (5.25,0) [fill=white] {};
     \node[vertex] (a3) at (7.25,0) [fill=white] {};
     \node[vertex] (b1) at (1.4,1.5) [fill=white] {};
     \node[vertex] (b2) at (3.75,1.5) [fill=white] {};
     \node[vertex] (b3) at (5.25,1.5) [fill=white] {};
     \node[vertex] (c1) at (0,3) [fill=white] {};
     \node[vertex] (c2) at (1.4,3) [fill=white] {};
     \node[vertex] (c3) at (2.9,3) [fill=white] {};
     \node[vertex] (c4) at (4.5,3) [fill=white] {};
     \node[vertex] (c5) at (6,3) [fill=white] {};
     \node[vertex] (c6) at (7.5,3) [fill=white] {};          
    \draw[very thick] (a1)--(b2)--(c3)  (a2)--(b3)--(c3) (a3)--(b3)
            (b1)--(c3)  (c1)--(c2)--(c3)--(c4)--(c5)--(c6);  
    \draw[thick] (-0.4,1) rectangle (3.3,3.6);    
    \draw[thick] (3.45,1) rectangle (5.55,2);    
    \draw[thick] (4.2,2.5) rectangle (7.8,3.6);                
          
    \node [above] at (3,3.1) {$u$};    
    \node [above] at (4.5,3.1) {$w$};        
    \node [above] at (-0.75,2.6) {$F_u$}; 
    \node [above] at (8.3,2.6) {$F_w$};     
    \node [above] at (5.75,1.3) {$A$};     
   \node [above] at (3.75,3.9) {$T$};         
    \end{scope}
\end{tikzpicture}   
%---------------------------------------------------------------
\begin{tikzpicture}
  [scale=0.7,inner sep=1mm, % this is the node radius
   vertex/.style={circle,thick,draw}, % this defines the default style for the vertex class
   thickedge/.style={line width=2pt}] % this defines the default style for the thick edges
    \begin{scope}[>=triangle 45]
    
     \node[vertex] (a1) at (3.75,0) [fill=white] {};
     \node[vertex] (a2) at (5.25,0) [fill=white] {};
     \node[vertex] (a3) at (7.25,0) [fill=white] {};
     \node[vertex] (b1) at (1.4,1.5) [fill=white] {};
     \node[vertex] (b2) at (3.75,1.5) [fill=white] {};
     \node[vertex] (b3) at (5.25,1.5) [fill=white] {};
     \node[vertex] (c1) at (0,3) [fill=white] {};
     \node[vertex] (c2) at (1.4,3) [fill=white] {};
     \node[vertex] (c3) at (2.9,3) [fill=white] {};
     \node[vertex] (c4) at (4.5,3) [fill=white] {};
     \node[vertex] (c5) at (6,3) [fill=white] {};
     \node[vertex] (c6) at (7.5,3) [fill=white] {};          
    \draw[very thick] (a1)--(b2)--(c4)  (a2)--(b3)--(c4) (a3)--(b3)
            (b1)--(c3)  (c1)--(c2)--(c3)--(c4)--(c5)--(c6);  
    \draw[thick] (-0.4,1) rectangle (3.3,3.6);    
    \draw[thick] (3.45,1) rectangle (5.55,2);    
    \draw[thick] (4.2,2.5) rectangle (7.8,3.6);                
          
    \node [above] at (3,3.1) {$u$};    
    \node [above] at (4.5,3.1) {$w$};        
    \node [above] at (-0.75,2.6) {$F_u$}; 
    \node [above] at (8.3,2.6) {$F_w$};     
    \node [above] at (5.75,1.3) {$A$};    
   \node [above] at (3.75,3.9) {$T'$};       
    \end{scope}
\end{tikzpicture}   
\caption{The trees $T$ and $T'$ in Lemma \ref{la:ironing-out}.}
\label{fig:ironing-1}
\end{center}
\end{figure}
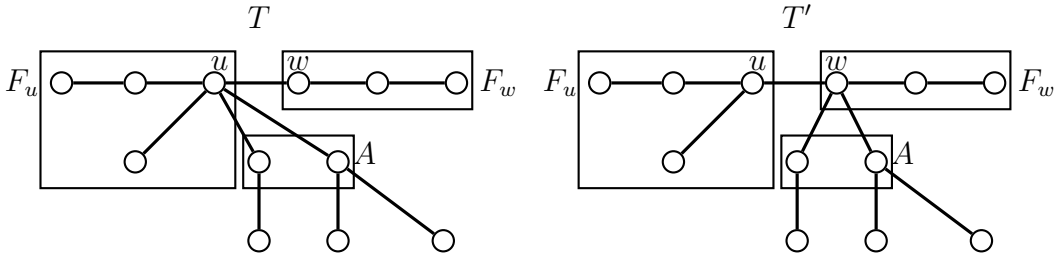

\begin{la} \label{la:ironing-out}
Given a tree $T$ with a weight function $c: V(T) \rightarrow \mathbb{N}_0$.
Let $u, w$ be two adjacent vertices of $T$. Let $A \subseteq N(u) -\{w\}$
be a nonempty set. Let $F$ be the component of 
$T-\{ua\ | \ a\in A\}$ containing $u$. Let $F_u$ and $F_w$ be the 
components of $F-uw$ containing $u$ and $w$,  
and $U$ and $W$ their vertex sets, respectively. 
Let $T'$ be the tree $T-\{ua \ | \ a\in A\} + \{wa \ | \ a\in A\}$.
If $c(U) > c(W)$, then
\[ SW_k(T,c) < SW_k(T',c). \]
\end{la}

{\bf Proof:}
For a sketch showing $T$ and $T'$ see Figure \ref{fig:ironing-1}. 
Let $V$ be the common vertex set of $T$ and $T'$, and let $S$ be a 
$k$-element subset of $V_c$. 
If $S^* \cap (U \cup W)= \emptyset$, then clearly $d_{T'}(S)=d_T(S)$. 
Similarly, if $S^*$ contains elements of both, $U$ 
and $W$, then  $d_{T'}(S)=d_T(S)$. Hence
$d_{T'}(S) \neq d_T(S)$ only if $S^* -(U\cup W) \neq \emptyset$ 
and exactly one of $S^*\cap U$ and $S^* \cap W$ is nonempty. Denoting
$V(T)-(U\cup W)$ by $X$ we obtain
\begin{eqnarray*}
SW_k(T',c) - SW_k(T,c) & = & \sum_{S\subseteq V_c, |S|=k} \big( d_{T'}(S)-d_T(S) \big) \\ 
   & \hspace*{-9em} = &  \hspace*{-8em}
   \sum_{S: S^*\cap X \neq \emptyset, S^*\cap U \neq \emptyset, S^*\cap W = \emptyset}
           \hspace*{-3em} \big(d_{T'}(S) - d_T(S)\big) + 
           \sum_{S: S^*\cap X \neq \emptyset, S^*\cap W \neq \emptyset, S^*\cap U = \emptyset}
           \hspace{-3em} \big(d_{T'}(S) - d_T(S)\big).           
\end{eqnarray*} 
If $S \subseteq V_c$ satisfies $S^* \cap X \neq \emptyset$, $S^* \cap U \neq \emptyset$ 
and $S^* \cap W = \emptyset$, then 
$d_{T'}(S) - d_T(S) = 1$. Similarly, 
if $S \subseteq V_c$ satisfies $S^* \cap X \neq \emptyset$, $S^* \cap W \neq \emptyset$ 
and $S^* \cap U = \emptyset$, then 
$d_{T'}(S) - d_T(S) = -1$. Hence
\begin{eqnarray*}
SW_k(T',c) - SW_k(T,c) & = &  \big| \{S\subseteq V_c \ | \ |S|=k, S^*\cap X \neq \emptyset, S^*\cap U \neq \emptyset, S^*\cap W = \emptyset \} \big|  \\
         & &   - 
           \big| \{S\subseteq V_c \ | \ |S|=k, S^*\cap X \neq \emptyset, S^*\cap W \neq \emptyset, S^*\cap U = \emptyset\} \big| \\
     & = & \sum_{i=1}^{k-1} \binom{c(U)}{i} \binom{c(X)}{k-i}   
          - \sum_{i=1}^{k-1} \binom{c(W)}{i} \binom{c(X)}{k-i}   \\
     & = & \sum_{i=1}^{k-1} \binom{c(X)}{k-i} 
          \Big[  \binom{c(U)}{i} - \binom{c(W)}{i} \Big] \\
     & > &  0,
\end{eqnarray*} 
where the last inequality holds since $c(U) > c(W)$. Hence the lemma follows. 
\hfill $\Box$ \\

The following lemma is central to our proof of the main results of
this paper. For trees it is a generalisation of inequality 
\eqref{eq:path-maximises-W}, which is the special case $k=2$ and $C=1$.

\begin{la}  \label{la:fundamental-lemma}
Let $T$ be a tree with vertex set $V$. Let 
$c:V \rightarrow \mathbb{N}_0$ be a weight function on the vertex
set of $T$ such that $c(v) \geq C$ for every vertex $v\in V$. 
Let $N$ be the total weight of the vertices of $T$. If $C\geq 1$, then
\[ SW_k(T,c) \leq \frac{k-1}{k+1} \frac{N+1}{C} \binom{N}{k} 
           +  \frac{C-1}{C} \binom{N}{k}, \] 
\end{la}

{\bf Proof:}
Let $N$ and $C$ be fixed. 
We may assume that among all trees and weight functions satisfying
the hypothesis of the lemma, $T$ and $c$ are
chosen such that $SW_k(T,c)$ is maximum. \\[1mm]
{\sc Claim 1:} $T$ is a path. \\
Suppose to the contrary that $T$ is not a path. Then $T$ contains
a vertex $u$ of degree at least three. Let $V_1, V_2, \ldots, V_r$ 
be the vertex sets of the components of $T-u$ and let $w_i$ the 
neighbour of $u$ in $V_i$ for $i=1,2,\ldots,r$. We may assume that
$c(V_1) \geq c(V_2) \geq \ldots \geq c(V_r)$. Let 
$A=\{w_1, w_2,\ldots,w_{r-2}\}$. Define the tree $T'$ by 
\[ T'= T-\{ua \ | \ a\in A\} + \{w_ra \ | \ a\in A\}. \]
Then the hypothesis of Lemma \ref{la:ironing-out}
is satisfied with vertices $u$ and $w_r$ corresponding to 
vertices $u$ and $w$ in Lemma \ref{la:ironing-out}, 
the subtree $T[V_{r-1} \cup V_r \cup \{u\}]$ corresponding
to $F$ in Lemma \ref{la:ironing-out}, and the sets
$\{u\} \cup V_{r-1}$ and $V_r$ corresponding to $U$ and $W$. 
We have 
$c(\{u\} \cup V_{r-1}) > c(V_{r-1}) \geq c(V_r)$.
Hence the hypothesis of Lemma \ref{la:ironing-out} is satisfied,
and if follows that
\[ SW_k(T',c) > SW_k(T,c). \]
which is a contradiction to the maximality of $SW_k(T,c)$. 
Hence $T$ is a tree with no vertex of degree greater than two, so $T$
is a path. This proves Claim 1. \\
From now on we assume that $T$ is a path with vertices $v_1, v_2,\ldots,v_n$
in this order. 
We now define an auxiliary graph $P$ on the vertex set $V_c$. 
Define the path $P$ by
\[   P = v_1^{1}, v_1^{2},\ldots,v_1^{c(v_1)},
         v_2^{1}, v_2^{2},\ldots,v_2^{c(v_2)},\ldots,
         v_n^{1}, v_n^{2},\ldots,v_n^{c(v_n)}. \]
{\sc Claim 2:} For any two vertices, 
$x,y$ of $P$ we have $d_T(x,y) \leq \frac{d_P(x,y)-1}{C} +1$. \\
Let $x=v_i^r$ and $y=v_j^s$. The statement clearly holds if $i=j$,
hence we assume, without loss of generality,  that $i<j$. 
Then, since $c(v_i), c(v_{i+1}),\ldots, c(v_j) \geq C$ we obtain  
\begin{eqnarray*} 
d_P(v_i^r, v_j^s) & \geq & d_P(v_i^{c(v_i)}, v_j^{1}) \\
             & = & c(v_{i+1}) + c(v_{i+2}) + \cdot + c(v_{j-1}) + 1 \\
             & \geq & C(j-i-1) + 1 \\ 
             & = & C(d_T(v_i,v_j)-1) +1,
\end{eqnarray*}
and by rearranging  we get Claim 2. \\[1mm]
{\sc Claim 3:} If $S \subseteq V_c$, $|S| \geq 2$, then  
$d_T(S) \leq \frac{d_P(S)-1}{C} +1$. \\
Since $T$ and $P$ are paths, we have
$d_P(S) = \max_{x,y \in S} d_P(x,y)$ and $d_T(S) = \max_{x,y \in S} d_T(x,y)$.
Hence, by Claim 2, 
\[ d_T(S) =  \max_{x,y \in S} d_T(x,y) \leq  \frac{\max_{x,y \in S} d_P(x,y)-1}{C} + 1
            =  \frac{d_P(S)-1}{C} +1, \]
which proves Claim 3. \\[1mm]
We now complete the proof. By Claim 3, 
\[ SW_k(T,c) = \sum_{S\subseteq V_c, |S|=k} d_T(S) 
     \leq \sum_{S\subseteq V_c, |S|=k} \Big( \frac{d_P(S)}{C} + \frac{C-1}{C} \Big)
     = \frac{1}{C} SW_k(P) + \frac{C-1}{C} \binom{N}{k}.  \] 
By \eqref{eq:path-maximises-W} we have $SW_k(P) = \frac{(k-1)(N+1)}{k+1}\binom{N}{k}$. 
Substituting this into the previous inequality we get
\[ SW_k(T,c) \leq \frac{1}{C}  \frac{(k-1)(N+1)}{k+1}\binom{N}{k}
           +  \frac{C-1}{C} \binom{N}{k}, \]
as desired. \hfill $\Box$

\section{A bound in terms of order and minimum degree}

In the proof of the following theorem we employ a refinement of a method that was 
developed in \cite{DanEnt2000} to prove an upper bound on the average distance. 
For a vertex $v$ of $G$ denote by $T(a)$ the subtree of $G$ whose vertex
set is $N[v]$, and whose edges are the edges joining $v$ to its neighbours
in $G$. 

\begin{theo}    \label{theo:bound-on-SW-for-min-degree}
Let $G$ be a connected graph of order $n$ with minimum degree $\delta$.
Then $G$ contains a spanning tree $T$ with
\[ SW_k(T) \leq \frac{k-1}{k+1} \frac{3(n+1)}{\delta+1} \binom{n}{k} 
           +  \Big(\frac{3\delta}{\delta+1}+2k\Big) \binom{n}{k}, \]
\end{theo}

{\bf Proof:} 
The strategy of the proof is as follows. We simultaneously construct a maximum
packing $A$ in $G$ and a subtree $T_0$ of $G$, which we extend to a spanning tree 
$T$ of $G$. We then show that the bound of the theorem holds for $T$.
 
To construct a packing of $G$ start by choosing a vertex $a_1$ and letting
$A_1=\{a_1\}$ and $T_1=T(a_1)$. Let $a_2$ be a vertex at distance exactly 
$3$ from $A_1$, if one exists. Then there exists an edge $e_2$ joining some vertex of 
$T_1$ to some vertex of $T(a_2)$. Let $A_2=A_1\cup \{a_2\}$ and 
let $T_2$ be the tree obtained from $T_1 \cup T(a_2)$ by adding 
the edge $e_2$. Let $a_3$ be a vertex at distance 
exactly $3$ from $A_2$, if one exists.  Then there exists an edge $e_3$ joining some 
vertex of $T_2$ to some vertex of $T(a_3)$. Let 
$A_3 = A_2\cup \{a_3\}$ and let $T_3$ be the graph obtained
from $T_2\cup T(a_3)$ by adding the edge $e_3$. 
Generally, given $A_i$ and $T_i$, we choose a vertex
$a_{i+1}$ at distance exactly $3$ from $A_i$, if one exists, let 
$e_{i+1}$ be an edge joining a vertex in $T_i$ to a vertex in
$T(a_{i+1})$, let
$A_{i+1}=A_i \cup \{a_{i+1}\}$,
and let $T_{i+1}$ be the tree obtained from the disjoint union of $T_i$ and  
$T(a_{i+1})$ by adding the edge $e_{i+1}$. 
Repeat this step until, after $s$ steps say,  all vertices are at distance 
at most $2$ from $A_s$. Let $T_0:=T_s$ and $A:=A_s$.  
Then $A$ is a maximal packing and every vertex of $G$ is within distance 
at most $2$ in $G$ from $A$, and thus adjacent to some vertex in $T_0$.
Joining every vertex not in $T_0$ to a neighbour in $T_0$ yields
a spanning tree $T$ of $G$. 

For every vertex $v$ of $T$ choose a nearest vertex  $a_v \in A$, 
and a shortest $(v,a_v)$-path $P_v$ in $T$. Then $P_v$ has at most
two edges. For $a\in A$ let $c(a)$ be the number of vertices $v$ 
of $G$ with $a_v=a$.
Since $A$ is a packing, for every $a\in A$ all vertices in $N[a]$ have
$a$ as their nearest vertex in $A$, and so  
\begin{equation} \label{eq:weight-at-least-delta+1}
 \textrm{$c(a) \geq \delta+1$ for all $a\in A$}, 
\end{equation} 
We now show that the difference $SW_k(T) -SW_k(T,c)$ is at most
$O(n^k)$. More specifically, we show that 
\begin{equation} \label{eq:difference-in-W-small}
SW_k(T) \leq  SW_k(T,c) + 2k \binom{n}{k}. 
\end{equation}
Let $c_1$ be the weight function that assigns a weight of $1$ to each vertex of $T$.
Clearly, $SW_k(T)=SW_k(T,c_1)$. Then $c$ is obtained by moving weight units from 
$v$ to $a_v$ for all vertices $v\in V$. Hence $V_c$ is obtained from $V_{c_1}$ by moving
the copy $v^1 \in V_{c_1}$ of vertex $v$ to vertex $a_v$ and making it a copy of $a_v$. 
Hence there is a bijection $f$ between the $k$-subsets of $V_{c_1}$ and 
the $k$-subsets of $V_c$, mapping every $k$-set $S_1 \subset V_{c_1}$
to a $k$-set $f(S) \subseteq V_{c}$ by replacing copies (with respect to $c_1$) 
of a vertex $v$ by copies of the vertex $a_v$ (with respect to $c$). 
If $T_S$ is a Steiner tree for a $k$-set $S \subseteq V_{c_1}$, then
by adding or deleting suitable edges that are in $\bigcup_{v \in S*} E(P_v)$, we
obtain a subtree of $T$ containing all vertices in $f(S)^*$. Since 
$|\bigcup_{v \in S*} E(P_v)| \leq 2k$, we conclude that
\[ d_T(S) \leq d_T(f(S)) + 2k.  \]
Summing over all $k$-subsets $S \subseteq V_{c_1}$ yields 
\eqref{eq:difference-in-W-small}.

We proceed to bound $SW_k(T,c)$. 
Since by the construction of $T_0$ and $T$ every vertex $a_i$, $i>1$, is 
at distance exactly $3$ 
in $T$ from some vertex $a_j$ with $j<i$,  it follows that 
in $T^3[A]$ there exists a path from $a_i$ to $a_1$ for every
$i>1$. Hence 
\begin{equation}  \label{eq:property-ii}
\textrm{ $T^3[A]$ is connected.}  
\end{equation} 
Let $H=T^3[A]$. Since $d_T(a_i, a_j) \leq 3d_H(a_i,a_j)$ for all
$a_i, a_j \in A$, and since the weight of $c$ is concentrated in
the vertices of $A$, we have
\begin{equation} \label{eq:Wiener-T-versus-H}
SW_k(T,c) \leq 3SW_k(H, c'), 
\end{equation}
where $c'$ is the restriction of $c$ to $A$. By \eqref{eq:weight-at-least-delta+1}, 
$c'(a_i) \geq \delta+1$ for all $a_i\in A$. Moreover, the total weight
of $c'$ is $n$, i.e., $N=n$. Hence, by Lemma \ref{la:fundamental-lemma}
\begin{equation}  \label{eq:obtained-from-fundamental-lemma}
SW_k(H,c') \leq \frac{k-1}{k+1} \frac{n+1}{\delta+1} \binom{n}{k} 
           +  \frac{\delta}{\delta+1} \binom{n}{k}. 
\end{equation} 
Applying \eqref{eq:difference-in-W-small}, \eqref{eq:Wiener-T-versus-H} and  
\eqref{eq:obtained-from-fundamental-lemma} 
we obtain
\[ SW_k(T) \leq \frac{k-1}{k+1} \frac{3(n+1)}{\delta+1} \binom{n}{k} 
           +  \Big(\frac{3\delta}{\delta+1}+2k\Big) \binom{n}{k}, \]
as desired. \hfill $\Box$

\begin{coro} \label{coro:bound-on-Wk-min-degree}
Let $G$ be a connected graph of order $n$ with minimum degree $\delta$.
Then 
\[ SW_k(G) \leq \frac{k-1}{k+1} \frac{3(n+1)}{\delta+1} \binom{n}{k} 
           +  \Big(\frac{3\delta}{\delta+1}+2k\Big) \binom{n}{k}, \]
           and thus
\[ \mu_k(T) \leq \frac{k-1}{k+1} \frac{3(n+1)}{\delta+1} 
           +  \frac{3\delta}{\delta+1}+2k. \]   
\end{coro}

\begin{exa} \label{exa-1} 
We now construct an example to show that the bound on the Steiner 
$k$-Wiener index in Corollary 
\ref{coro:bound-on-Wk-min-degree} is best possible apart from 
a term $O(n^k)$, and that the bound on the average Steiner $k$-distance
is best possible apart from an additive constant. We only
construct examples for the case that $\delta+1$ is a multiple
of $3$, but it is not difficult to modify this
construction for all values of $\delta$. For $d \in \mathbb{N}$ 
define the graph $G_{d,\delta}$ by 
\[ G_{d,\delta} = K_{\delta} + K_{(\delta+1)/3} + K_{(\delta+1)/3} 
           + \cdots + K_{(\delta+1)/3} + K_{(\delta+1)/3} 
            + K_{\delta}, \]
where the term $K_{(\delta+1)/3}$ appears $d-1$ times. 
Clearly, $n(G_{d,\delta}) = \frac{d+5}{3}(\delta+1)-2$, and 
${\rm diam}(G_{d,\delta}) = d$. Hence for large $d$ and constant $\delta$ we have 
$n=d\frac{\delta+1}{3} +O(1)$. We now bound the Steiner $k$-Wiener
index from below. For $i=1,2,\ldots,d-1$ let  $V_i$ be the set of
vertices of the $i$-th copy of $K_{(\delta+1)/3}$, and let
$V_0$ and $V_d$ be subsets of cardinality $(\delta+1)/3$ of
the first and last, respectively, copy of $K_{\delta}$. Let  
${\cal S}$ be the set of all $k$-element sets of vertices of $G_{d,\delta}$
that  are contained in $\bigcup_{i=0}^d V_i$, and that 
have no two vertices in the same $V_i$.  
Let $F$ be the path of order $d+1$ with vertices
$u_0,u_1,\ldots,u_d$. We define a mapping $f$ that maps every set in ${\cal S}$
to a $k$-set of vertices of $F$. For $S \in {\cal S}$ let 
$f(S)$ be the subset of $V(F)$ containing those $u_i$ for 
which $S$ contains a vertex in $V_i$. It is clear that
$d_{G_{d,\delta}}(S) = d_{F}(f(S))$. Since every $k$-set
of vertices of $F$ is the image under $f$ of exactly
$(\frac{\delta+1}{3})^{k}$ sets in ${\cal S}$, we have 
\[ SW_k(G_{d,\delta}) \geq \sum_{S \in {\cal S}} d_{G_{d,\delta}}(S) 
            = \sum_{S \in {\cal S}} d_{P_{d+1}}(f(S)) 
            =  \sum_{S \subseteq V(F), |S|=k} (\frac{\delta+1}{3})^{k} d_F(S).
           \]
Hence, by Theorem \ref{theo:bounds-on-Wk-for-trees}, 
\[  
SW_k(G_{d,\delta})  \geq  (\frac{\delta+1}{3})^{k}\, SW_k(F) \\
               =   (\frac{\delta+1}{3})^{k}           
                 \frac{(k-1)(d+2)}{k+1} \binom{d+1}{k}.
\]
Now for constant $k$ and $\delta$ and large $n$ and $d$ we get  
$n=d\frac{\delta+1}{3} +O(1)$ and thus $d=\frac{3n}{\delta+1}+O(1)$. 
Hence $\binom{d+1}{k} = (\frac{3}{\delta+1})^k \binom{n}{k} + O(n^{k-1})$, 
and so
\begin{eqnarray*} 
SW_k(G_{d,\delta})  & \geq &  (\frac{\delta+1}{3})^{k}           
     \frac{(k-1)(d+2)}{k+1}  \Big[ (\frac{3}{\delta+1})^k   \binom{n}{k} +O(n^{k-1}) \Big]  \\
     & = & \frac{k-1}{k+1} \frac{3n}{\delta+1} \binom{n}{k} +O(n^k). 
\end{eqnarray*} 
Dividing by $\binom{n}{k}$ we get 
\[ \mu_k(G_{d,\delta}) = \frac{k-1}{k+1} \frac{3n}{\delta+1} +O(1), \]
as desired.
\end{exa}

\section{An improved bound for triangle-free graphs} 

Our main aim in this section is to improve the bound in 
Theorem \ref{theo:bound-on-SW-for-min-degree} for triangle-free graphs. 
The basic idea of the proof of the improved bound is similar to  
Theorem \ref{theo:bound-on-SW-for-min-degree}, but some additional arguments are needed. 

For an edge $e=uv$ of a triangle-free graph $G$ denote by $T(e)$ the subtree of 
$G$ whose vertex set is $N(u) \cup N(v)$, and whose edges are the edges joining 
$u$ or $v$ to its neighbours in $G$. The distance $d(e_1,e_2)$ between two
edges of $G$ is the minimum of the four distances between a vertex incident 
with $e_1$ and a vertex incident with $e_2$. The distance between an edge
$e$ and a set $E_1$ of edges is the minimum of the distances between 
$e$ and the edges in $E_1$.

\begin{la} \label{la:line-graph-distance}
Let $T$ be a tree and $L$ the line graph of $T$. Let $S_V$ be a set of vertices
of $T$, and $S_E$ a set of edges of $T$ such that $S_V \subseteq V(S_E)$. Then
\[ d_T(S_V) \leq d_L(S_E) +1. \]
\end{la}

{\bf Proof:} Let $T_L(S_E)$ be a Steiner tree for $S_E$ in $L$ 
and let $U$ be its vertex set. Then $U$ is a 
set of $d_L(S_E)+1$ vertices of $L$. Since 
$U\subseteq E(T)$ and since $U$ induces a connected graph in $L$, the subgraph 
$T[U]$ of $T$ induced by the set of edges $U$ is also connected. Since
$T[U]$ contains all vertices of $S_V$, we have $d_T(S_V) \leq |U| = d_L(S_E)+1$,
as desired. 
\hfill $\Box$

\begin{theo}    \label{theo:triangle-free-bound-on-SW-for-min-degree}
Let $G$ be a connected, triangle-free graph of order $n$ with minimum degree $\delta$.
Then $G$ contains a spanning tree $T$ with
\[ SW_k(T) \leq \frac{k-1}{k+1} \frac{2(n+1)}{\delta} \binom{n}{k} 
           +  \Big(\frac{4\delta-2}{\delta}+3k+1\Big) \binom{n}{k}. \]
\end{theo}

{\bf Proof:} 
To construct a matching of $G$ start by choosing an edge $b_1$ and letting
$M_1=\{b_1\}$ and $T_1=T(b_1)$. Let $b_2$ be an edge at distance exactly 
$3$ from $M_1$, if one exists. Then there exists an edge $e_2$ joining some vertex of 
$T_1$ to some vertex of $T(b_2)$. Let $M_2=M_1\cup \{b_2\}$ and 
let $T_2$ be the tree obtained from $T_1 \cup T(b_2)$ by adding 
the edge $e_2$. Let $b_3$ be an edge at distance exactly $3$ from 
$M_2$, if one exists.  Then there exists an edge $e_3$ joining some 
vertex of $T_2$ to some vertex of $T(b_3)$. Let 
$M_3 = M_2\cup \{b_3\}$ and let $T_3$ be the graph obtained
from $T_2\cup T(b_3)$ by adding the edge $e_3$. 
Generally, given $M_i$ and $T_i$, we choose an edge 
$b_{i+1}$ at distance exactly $3$ from $M_i$, if one exists, let 
$e_{i+1}$ be an edge joining a vertex in $T_i$ to a vertex in
$T(b_{i+1})$, let
$M_{i+1}=M_i \cup \{b_{i+1}\}$,
and let $T_{i+1}$ be the tree obtained from $T_i\cup T(b_{i+1})$
by adding the edge $e_{i+1}$. 
Repeat this step until, after $s$ steps say,  all edges are at distance 
at most $2$ from $M_s$. Let $T_0:=T_s$ and $M:=M_s$.  
Then $M$ is a matching and every edge of $G$ is within distance 
at most $2$ in $G$ from $M$, and so every vertex of $G$ is at
distance at most three from $V(M)$.
Joining every vertex not in $T_0$ to a neighbour that is closer to 
$T_0$ or in $T_0$ yields a spanning tree $T$ of $G$ that
preserves the distances from all vertices to $V(M)$. 

We now show that $SW_k(T)$ is bounded as claimed. 
For every vertex $v$ of $T$ choose a nearest vertex  $a_v \in V(M)$, 
and a shortest $(v,a_v)$-path $P_v$ in $T$. Then $P_v$ has at most
three edges. For $a\in V(M)$ let $c(a)$ 
be the number of vertices $v$ of $G$ with $a_v=a$.
Since $M$ is a matching and since $G$ is triangle-free, for every 
$a\in V(M)$ all vertices in $N_G[a]$ except the matching partner of $a$ 
have $a$ as their nearest vertex in $M$, and so  
\begin{equation} \label{eq:weight-at-least-delta-triangle-free}
 \textrm{$c(a) \geq \delta$ for all $a\in M$}, 
\end{equation} 
Making use of the fact that every vertex of $T$ is within distance
three of some vertex in $V(M)$, we show 
as in the proof of Theorem \ref{theo:bound-on-SW-for-min-degree} 
(see equation \eqref{eq:difference-in-W-small} there) we show 
\begin{equation} \label{eq:difference-in-W-small-triangle-free}
SW_k(T) \leq  SW_k(T,c) + 3k \binom{n}{k}. 
\end{equation}
Let $L(T)$ be the line graph of $T$. 
Define a weight function $c_2$ on the vertices of $L(T)$, i.e.,
the edges of $T$, by
\[ c_2(uv) = \left\{ \begin{array}{cc}
        c(u) + c(v) & \textrm{if $uv \in M$}, \\
        0 & \textrm{if $uv \notin M$.}
            \end{array} \right. \] 
Clearly, $c(V(T))=c_2(M)=n$. We now show that 
\begin{equation}  \label{eq:WT-vs-WL}
SW_k(T, c) \leq SW_k(L,c_2) + \binom{n}{k}. 
\end{equation} 
Define a bijection $f:V(T)_c \rightarrow V(L)_{c_2}$  that maps, 
for every edge $uv\in M$, 
$\{u^1, u^2\ldots, u^{c(u)}\} \cup \{v^1, v^2\ldots, v^{c(u)}\}$ 
to $\{uv^1, uv^2\ldots, uv^{c_2(uv)}\}$. If $S \subseteq V_c$ 
then $f(S) \subseteq E(T)_{c_2}$ and clearly 
$S^* \subseteq V(f(S)^*)$. 
Hence, by Lemmal \ref{la:line-graph-distance} this implies 
\[ d_T(S^*) \leq d_L(f(S)^*) +1, \]
and so   
\[ d_T(S) \leq d_L(f(S)) +1, \]
Summation over all $k$-element subsets $S$ of $V_c$ yields
\begin{eqnarray*}
\sum_{S \subseteq V_c, |S|=k} d_T(S) 
    & \leq & \sum_{S \subseteq V_c, |S|=k} \big(d_L(f(S))+1\big) \\
      & = &  \sum_{S \subseteq E(T)_{c_2}, |S|=k} \big(d_L(f(S))+1\big) \\
      & = & SW_k(L,c_2) + \binom{n}{k},
\end{eqnarray*}
which is \eqref{eq:WT-vs-WL}.

By the construction of $T_0$ and $T$ every edge $b_i$, $i>1$, is 
at distance exactly three 
in $T$ from some edge $b_j$ with $j<i$. It follows that in the line graph $L$
every vertex $b_i$ of $L$ with $i>1$ is at distance exactly four  from some 
vertex $b_j$ with $j<i$. Therefore, 
\begin{equation}  \label{eq:property-ii-triangle-free}
\textrm{ $L^4[M]$ is connected.}  
\end{equation} 
Let $H=L^4[M]$. Since $d_L(b_i, b_j) \leq 4d_H(b_i,b_j)$ for all
$b_i, b_j \in M$, and since the weight of $c_2$ is concentrated in
$M$, we have
\begin{equation} \label{eq:Wiener-T-versus-H-triangle-free}
SW_k(L,c_2) \leq 4\, SW_k(H, c'), 
\end{equation}
where $c'$ is the restriction of $c$ to $M$. By 
\eqref{eq:weight-at-least-delta-triangle-free}, 
$c'(b_i) \geq 2\delta$ for all $b_i\in M$. Moreover, the total weight
of $c'$ is $n$, i.e., $N=n$. Hence, by Lemma \ref{la:fundamental-lemma}
\begin{equation}   \label{eq:obtained-from-fundamental-lemma-triangle-free}
SW_k(H,c') \leq \frac{k-1}{k+1} \frac{n+1}{2\delta} \binom{n}{k} 
           +  \frac{2\delta-1}{2\delta} \binom{n}{k}. 
\end{equation}           
Applying \eqref{eq:difference-in-W-small-triangle-free}, \eqref{eq:Wiener-T-versus-H-triangle-free} 
and \eqref{eq:obtained-from-fundamental-lemma-triangle-free} we obtain
\[ SW_k(T) \leq \frac{k-1}{k+1} \frac{2(n+1)}{\delta} \binom{n}{k} 
           +  \Big(\frac{4\delta-2}{\delta}+3k+1\Big) \binom{n}{k}, \]
as desired. \hfill $\Box$

\begin{coro} \label{coro:bound-on-Wk-min-degree-triangle-free}
Let $G$ be a connected graph of order $n$ with minimum degree $\delta$.
Then 
\[ SW_k(G) \leq \frac{k-1}{k+1} \frac{2(n+1)}{\delta} \binom{n}{k} 
           +  \Big(\frac{4\delta-2}{\delta}+3k+1\Big) \binom{n}{k}, \]
           and thus
\[ \mu_k(G) \leq \frac{k-1}{k+1} \frac{2(n+1)}{\delta} 
           +  \frac{4\delta-2}{\delta}+3k+1. \]
\end{coro}

\begin{exa} 
The following example shows that the bound on the Steiner 
$k$-Wiener index in Corollary 
\ref{coro:bound-on-Wk-min-degree-triangle-free} is best possible apart from 
a term $O(n^k)$, and that the bound on the average Steiner $k$-distance
is best possible apart from an additive constant. We only
construct examples for the case that $\delta$ is even, 
but as in Example \ref{exa-1} it is not difficult to modify this
construction for odd values of $\delta$. For $d \in \mathbb{N}$ 
define the graph $H_{d,\delta}$ by 
\[ H_{d,\delta} = \delta K_1 + \delta K_1 + K_{\delta/2} + K_{\delta/2} 
     + \cdots + K_{\delta/2} + K_{\delta/2}  + \delta K_1 + \delta K_1, \]
where the term $K_{\delta/2}$ appears $d-3$ times. 
Then calculations similar to those in Example \ref{exa-1} show that for
constant $\delta$ and $k$ and large $n$ and $d$ we have 
\[ SW_k(H_{d,\delta}) \leq \frac{k-1}{k+1} \frac{2(n+1)}{\delta} \binom{n}{k} 
           +  O(n^k), \]
           and thus
\[ \mu_k(H_{d,\delta}) \leq \frac{k-1}{k+1} \frac{2n}{\delta} 
           +  O(1). \]
\end{exa}

\end{document}